%% file: simplearchimedean
\input amstex

\input generic_macros

\input papermacros_entirefunctions

\def\Z{\text{\bf Z}}


\def\Z{\text{\bf Z}}

\long\def\subsubtitle #1.{\vskip 3pt \noindent {\it #1 }}
\def\Aff{\text{Aff\,}}

\long\def\subtitle #1\par#2{\vskip 3pt \noindent {\it #1}\par \noindent
#2}

\def\ideal(#1){\(#1\)}
\def\oid(#1){< #1 >}

\Title Simple archimedean dimension groups%
\plainfootnote{\rm${}^0$}{\rm Working document.}

\comment
\Abstract The main result is that dimension groups (not countable) that are also real ordered vector spaces can be obtained as  direct limits (over  directed sets) of simplicial real vector spaces (finite dimensional vector spaces with the coordinatewise ordering), but the directed set is not as interesting as one would like---e.g., it is not clear whether a countable dimensional example can be represented as such a direct limit over the positive integers. It turns out this is the case when the group is additionally simple, and it is shown that the latter have an ordered tensor product decomposition. 
\endcomment

\vskip4pt 
\noindent {\it David Handelman}\plainfootnote{$^1$}{Supported in part
by a Discovery Grant from NSERC.}\vskip 4pt 

\noindent
Let $(G,u)$ be an unperforated partially ordered abelian group with an order unit, $u$. It is {\it simple\/} if every nonzero positive element is an order unit. By [G, ], it is {\it archimedean\/}%
\plainfootnote{$^2$}{This uses the classical definition of archimedean: for elements $g$ and $h$ of $G$, $ng \leq h$ for all positive integers $n$ implies  $-g \geq 0$. In the presence of an order unit, this  is equivalent to the trace-determining condition in the text, and in particular shows that the map $G \to \Aff S(G,u)$ is an embedding (an order-embedding, in fact) when $G$ is archimedean. A much weaker definition---irrelevant here---is used by a large group of workers in real algebraic geometry, causing confusion. }
if and only if for $g$ in $G$, $\tau (g) \geq 0$ for all pure traces $\tau$ implies $g \in G^+$. Among unperforated groups, archimedean and simple represent properties that are maximal and minimal, respectively---the former are those for which the weakest necessary condition (that the values at traces be nonnegative) implies positivity, while the latter are those for which the strongest sufficient condition (that positivity implies strict positivity on traces) is implied by positivity. So it is a little difficult to construct dimension groups that are both simple and archimedean, aside from subgroups of the reals. Question 4 of [G] asks whether every Choquet simplex can be the trace space a of a simple archimedean dimension group; we show this is the case for metrizable simplices, based on an interesting construction over the  interval.

The normalized trace space of $(G,u)$ is denoted $S(G,u)$; the latter's extremal boundary is denoted $\partial_e S(G,u)$. We use $\,\hat{}\,$ to denote the natural map $G\to \Aff S(G,u)$, given by $\hat g (\tau) = \tau (g)$. 

An alternative formulation of simple and archimedean (in the presence of an order unit) is the following. Suppose $G$ admits an unperforated partial ordering \st the map $G \to \Aff S(G,u)$  is an embedding; then the ordering is unique (the smallest ordering is strict ordering, the largest ordering is the pointwise one).

The following is practically tautological.

\Lem Lemma 1.  Let $(G,u) $ be an unperforated  partially ordered abelian group with order unit. Then $G$ is simple and archimedean if and only for all $g \in G\setminus\brcs{0}$, $\inf_{\tau \in \partial_e S(G,u)} \tau (g) \neq 0$.

\Pf A standard facial argument shows that  $\inf_{\tau \in \partial_e S(G,u)} \tau (g) > 0$ implies $\inf_{\tau \in  S(G,u)} \tau (g) > 0$ (even though $\partial_e S(G,u)$ need not be compact)---the condition implies that  $\tau (g) \geq 0$ for all traces $\tau$;  set $F = \Set{\tau \in S(G,u)}{\tau(g) =0 }$---it is easy to see that if nonempty, this is a closed face, hence has extreme points, which (since $F$ is a face) are extreme \wrt $S(G,u)$. 

Assume $G$ is archimedean and simple; the former says that  $\inf_{\tau \in \partial_e S(G,u)} \tau (g) = 0$ entails $g \in G^+$; simplicity implies  $g$ would be an order unit, hence the infimum would be strictly greater than zero. 

Conversely, suppose nonzero $g$ satisfies $\tau (g) \geq 0$ for all pure $\tau$. By hypothesis (the infimum is not zero, hence must be strictly positive), $\tau (g)>  0$ for all pure $\tau$, hence (by the first paragraph), $\tau (g) > 0$ for all traces, whence $g$ is an order unit and thus in the positive cone, so $G$ is simple and archimedean (simultaneously).
\qed

Note that the criterion refers to {\it all\/} nonzero $g$, not just those in $G^+$ (which would characterize simplicity). In particular,  a simple archimedean group which is also an ordered real vector space must be the reals with the usual ordering (pick any nonzero $g$ in $G^+$, and let $\alpha = \inf \tau(g)$; then $g-\alpha u$ vanishes at an extreme point---again, using the facial argument---and is nonnegative everywhere; archimedeanness entails $g \geq \alpha u$, and simplicity then forces $g = \alpha u$).

An extreme version is the following. An unperforated partially ordered abelian group with order unit, $(G,u)$, is {\it extremely simple\/} if for all $g \in G \setminus \brcs{0}$, and all pure traces $\tau$, $\tau(g) \neq 0$. 

To see that this implies the criterion of Lemma 1, we use the facial remark in the first paragraph of the argument---if $G$ is extremely simple and $\tau (g) \geq 0$ for all pure $\tau$, the hypothesis ensures that $\tau(g) > 0$ for all pure $\tau$, and the facial argument entails $g$ is an order unit, so the infimum is strictly positive. 

 Obvious  examples are subgroups of the reals with the relative ordering. There are others.

\Lem .Examples of extremely simple dimension groups with modestly interesting trace spaces. 

\noindent 1 Suppose $\alpha_i$ ($i=1,2,\dots,n$) are real numbers \st each of $\brcs{1, \alpha_i}$ is linearly independent over the rationals (that is, none of $\alpha_i$ is rational). Let $G$ be the subgroup of $\R^n$ spanned (as a group) by the standard basis vectors $\brcs{e_i}$ together with the element $E = \sum \alpha_j e_j$. Then $G$ is a free abelian group of rank $n+1$; equipped with the relative order inherited from $\R^n$   (i.e., $G^+ = (\R^n)^+ \cap G$, where $\R^n$ has the usual coordinatewise ordering), $G$ is unperforated and it is easy to check that all the pure traces on $G$ are given by  the $n$ coordinate evaluations, and the linear independence hypotheses ensure $G$ is extremely simple. In this case, the pure trace space consists of $n$ points.

If we additionally insist that the set $\brcs{1, \alpha_1, \alpha_2, \dots, \alpha_n}$ be linearly independent over the rationals, then as is well known, $G$ is a dense subgroup of $\R^n$, and being simple, is thus a dimension group. 

\vskip 4pt\noindent 2 Suppose that $K$ is an algebraic extension field of $\Q$, and in addition $K$ is formally real, i.e., if not all $k_i$ are zero, then $\sum k_i^2 \neq 0$ (or what amounts to the same thing, $K$ admits a real embedding). Impose on $K$ the sums of squares ordering (that is, $K^+$ consists of the set of sums of squares). As $K$ is formally real, this is a proper cone, and from algebraicity, it follows and is easy to check that $1$ is an order unit for $K$; moreover, inverses of positive elements are positive, and multiplication preserves the positive cone. It is known that $K$ is a dimension group [H; old paper], i.e., satisfies interpolation (this is true for any formally real field, not just algebraic extensions of the rationals). For a partially ordered ring with $1$ as order unit, the pure traces are multiplicative, in particular are ring homomorphisms. Since $K$ is a field, none of these have  nontrivial kernel, verifying extreme simplicity. In particular, $K$ is a simple archimedean dimension group, which is also an ordered ring.

The pure trace space can be interesting. Let  $K$ be $\Q_{\R}$,  the subfield of the reals consisting of all elements algebraic over $\Q$ ($\Q_{\R}$ is of index two in the algebraic closure of $\Q$). Then the pure traces can be identified with the Galois automorphisms of $\Q_{\R}$, and in particular, the pure trace space is the infinite, nonatomic,  and totally disconnected separable set (sometimes called the (or a) Cantor set, although {\it von Neumann compactum\/} would be at least as  appropriate since he proved its uniqueness). \qed

Extreme simplicity is drastic, as evidenced by the trivial Lemma 2.

\Lem Lemma 2. Suppose that $(G,u)$ is a dimension group with order unit, and there exists a connnected subset $U$ of $\partial_e S(G,u)$. Then there exists  $g$ in $G$ and $s$ in $U$ \st $\hat g$ is not constant on $U$, and $s(g) = 0$.

\Pf  Given distinct $v$ and $w$ in $U$, there exists $h$ in $G$ \st $\hat h (v) \neq \hat h (w)$. Since $\hat h|U$ is continuous and $U$ is connected, the range of $\hat h|U$ is a nontrivial interval. Let $q = a/b$ be a rational number (with $a$ an integer and $b$ a positive integer) in the interval. There exists $s$ in $U$ \st $\hat h (s) = q$. Set $g = bh - au$; since $\hat h$ is nonconstant on $U$ and $\hat u $ is the constant function $1$, $\hat g$ is not constant on $U$.\qed

In particular, if $(R,1)$ is an unperforated partially ordered ring with $1$ as order unit, and $R$ is extremely simple, then the pure trace space (known to be compact, since the pure traces are exactly the multiplicative ones) must be totally disconnected, as in Example 2 above. We simply note that if $\partial_e S(R,1)$ contained a connected subset, then by Lemma 2, there would be an element $r$ together with a pure trace $x$ \st $x(r) = 0$, and moreover archimedeanness entails $\hat r \neq 0$, so $r^2$ is not zero. Then $r^2$ is nonnegative at every pure trace (since all such are multiplicative), hence by archimedeanness, $r^2 $ would belong to $R^+$, and of course, $x(r^2) = 0$, contradicting Lemma 1. 

This seems about the end of the road for extremely simple (dimension) groups.

\Lem Example 3. Simple archimedean dimension groups with  the unit interval as pure trace space. 

\Pf Recall from Example 2 above, the maximal algebraic (over the rationals) subfield of the reals, $\Q_{\R}$. Let $\brcs{\alpha_i}_{i \in \N}$ be a countably infinite set of real numbers \st the enlarged set $\brcs{1} \cup \brcs{\alpha_i}_{i \in \N}$ is linearly  independent over $\Q_{\R}$---e.g., if $t$ is a transcendental number, we could take $\alpha_i = t^i$, or we could simply insist that $\brcs{\alpha_i}$ be {\it algebraically\/} independent over $\Q$. 

Inside the real polynomial algebra $\R[x]$, define the elements, $e_0 = 1$, $e_i = x^i - \alpha_i$ ($i=1,2,3, \dots$), and define $G \subset \R[x]$ to be the rational span of $\brcs{e_j}_{j \geq 0}$.

Equip $\R[x] $ with the strict ordering as a subgroup of $C([0,1],\R)$, so that $\R[x]$ is a simple dimension group (since the image is dense), and put the relative ordering on $G$. Automatically, $G$ is simple. Next, $G$ is a rational vector space, so its closure \wrt the supremum norm---equivalently the norm on $\R[x]$ with the strict ordering---is a real vector space, and thus each  $x^i$ is contained in the closure of $G$. Hence $G$ is dense in $C([0,1],\R)$ (since $\R[x]$ is), and in particular, its pure trace space is the same as that of $\R[x]$, namely $[0,1]$, and moreover, $G$ is a dimension group.

 We show that $G$ is archimedean by verifying the condition of Lemma 1.

Pick (to begin with) an arbitrary nonzero element of $G$, $g = q_0 + \sum_1^n  (x^i - \alpha_i) q_i$, where $q_i$ are rationals; since $g$ is not zero, not all the rational coefficients are zero. If $\alpha$ is a real number that is algebraic over the rationals, then $g (\alpha ) = 0$ entails $q_0 + \sum_1^n q_i( \alpha^i - \alpha_i) =0$, yielding the equation, $q_0 + \sum_{i=1} q_i \alpha^i = \sum q_i \alpha_i$; the left side is algebraic (over the rationals), hence belongs to $\Q_{\R}$, while the right side is a rational- (hence a $\Q_{\R}$-) linear combination of $\brcs{\alpha_i}$. By our assumption, both sides must be zero, which forces $q_1 = q_2 = \dots =0$, and this in turn forces $q_0 = 0$, a contradiction. The conclusion is that if $g$ is   an element of $G$ and not a constant, then it cannot have any zeros at algebraic points.

Now suppose that nonzero $g$ in $G$ has minimum $0$ as a function on the unit interval. Then  $\inf_{}g(\alpha) = 0$ for some $\alpha$ in the unit interval. By the preceding paragraph, $\alpha$ is not algebraic, so in particular, cannot be zero or one; thus it must be an interior point, and since $g$ is nonnegative, $\alpha$ is the location of  a minimum of $g$ (as a continuous function on $[0,1]$). Since $g$ is a polynomial and $\alpha$ is an interior point of the interval, we must have $g'(\alpha) = 0$. But this entails
$\sum q_i i \alpha^{i-1} = 0$, which in turn entails either that $\alpha$ is algebraic, or that $q_1 = q_2 = \dots = q_n = 0$, hence $q_0 = 0$; either way, we reach a contradiction. Thus $g$ in $G\setminus\brcs{0}$ with $g|[0,1] \geq 0$ forces $g$ to have no zeros in $[0,1]$.
This verifies the criterion of Lemma 1. 

 Thus $G$ is a simple archimedean dimension group whose pure trace space is the unit interval. 

A  sensitivity phenomenon occurs if we change the endpoints of the interval from $[0,1]$ to $[a,b]$. If both $a$ and $b$ are algebraic, then the same argument applies (since the putative functions cannot vanish at either endpoint, hence any zeros must be in the interior, whence the derivative argument works). On the other hand, if either $a$ or $b$ is of the form $\alpha_k^{1/k}$ (where this makes sense, e.g., if $k$ is even, then $\alpha_k$ must be positive) for some $k$, then one of  $\pm e_k =\pm(x^k - \alpha_k)$ (an element of $G$) will vanish at one point of the interval while being strictly positive on the rest of it---in particular, the so-constructed simple dimension group $G$ will {\it not\/} be archimedean.
\qed

In this example, we can consider the ordered tensor product, $H := G \otimes_{\Z} \R$; as vector spaces, the obvious map $H \to \R[x]$ is an isomorphism inducing an affine homeomorphism on their respective trace spaces; it is also positive. Since $G$ and $\R$ are simple dimension groups, so is $H$, and it follows that the map is an order-isomorphism (of ordered vector spaces).

However, we could just as well have imposed a different ordering on $\R[x]$ which yields exactly the same pure traces---for example, the positive cone generated additively and multiplicatively by $\brcs{\R^+, x, 1-x}$ (Renault's example). This is a non-simple dimension group with pure trace space $[0,1]$. The inclusion $G \subset \R[x]$ yields the same ordering on $G$ (that is, with this latter ordering on $\R[x]$, the relative ordering on $G$ is the same as the original), so that in this case the natural map $H \to \R[x]$ is not an order isomorphism (the left is simple, the right isn't), although it is a vector space isomorphism which is also positive.

Another candidate for the ordering on $\R[x]$ is pointwise on $[0,1]$, that is, make $\R[x]$ itself archimedean; the same remarks apply, except I cannot see whether it is   a dimension group.

The idea underlying Example 3 yields a complete answer (at least in the metrizable case---without metrizability, there probably is an argument, but it looks like a lot of effort) to Goodearl's question. 

\Lem Example $\infty$. For every metrizable Choquet simplex $K$, there exists a countable simple archimedean dimension group whose trace space is  (affinely homeomorphic to)  $K$. 

\Pf Let $K$ be a metrizable Choquet simplex, and let $u_0 = 1$, $u_1$, \dots  be a countable set of elements of $A := \Aff K$ (where $1$ simply means the constant function) \st $\brcs{u_i}$ is linearly independent over the reals, and its  real span    is dense in $A$. We will construct out of this a rational vector subspace of $A$ (parallel to the development of the polynomial example), $G$, satisfying the criterion of Lemma 1. 

Let $s_- (g) = \inf_{\tau \in \partial_e K} g(\tau)$ and $s_+ (g) = \sup_{\tau \in \partial_e K} g(\tau)$. The facial argument yields that the values are actually attained on $\partial_e K$, i.e., there exist $\tau_+$ and $\tau_-$ in $\partial_e K$ \st both equations $\tau_{\pm} (g) = s_{\pm} (g)$ hold. For a subgroup $J$ of $A$, denote by $s(J)$ the subgroup of the reals generated by  set of values of $s_- (g)$ as $g$ varies over $J$; since $J$ is closed under multiplication by $-1$, this is the same as the group generated by the set of values of $s_+(g)$, which explains the lack of sign in $s(J)$. Obviously $s(J)$ is a subgroup of the reals; if $J$ is a rational vector space, so is $s(J)$, and if $J$ is countable, so is $s(J)$. It is not clear that the set of values of the $s_-(g)$ (running over $J$) {\it is\/} a group, but fortunately all that matters for this example is cardinality of the group it generates.  

We proceed to define $v_i$ inductively, so that if $H_i $ is the rational vector space span of $\brcs{v_0, v_1, \dots, v_i}$, then $\R H_i$ is the real span of $\brcs{u_0, u_1, \dots, u_i}$ for all $i$, and various other properties. Let $v_0 = u_0 =1$; suppose we have $v_0,\dots,v_{k-1}$ with the following properties:

\item{(a)} there exist nonzero real numbers $\lambda_i $ \st for all $1 \leq i\leq k-1$, $v_i = u_i - \lambda_i$ (this notation identifies $\lambda_i$ with the corresponding constant function; to be pedantic, $v_i = u_i - \lambda_i {\pmb 1}$ where $\pmb 1$ is the constant function with value $1$);

\item{(b)} On defining $H_i $ as above, for $1 \leq i \leq k-1$, for all $g \in H_i \setminus H_{i-1}$,  each of $s_{\pm} (g)$
 is a nonzero rational multiple of $ \lambda_i$ modulo $s(H_{i-1} + \Q u_i)$, and moreover, $\lambda_i \Q \cap s(H_{i-1} + \Q u_i) = \brcs{0}$.

\noindent We will show the inductive process can be continued.

Consider the  countable rational vector subspace of the reals, $s(H_{k-1} + \Q u_k)$; we may thus select nonzero real $\lambda_k$ \st $\lambda_k \Q \cap s(H_{k-1} + \Q u_k) = \brcs{0}$, and set $v_k = u_k - \lambda_k$. Obviously the real span of $\brcs{u_i}_{i=0}^k$ coincides with the real span of $\brcs{v_i}_{i=0}^k$. We observe that since the set $\brcs{u_j}$ is linearly independent over the reals, so is the set $\brcs{v_0, v_1,\dots,v_k}$, and therefore it is linearly independent  over the rationals. Then define $H_k = H_{k-1} + v_k \Q$, and select $g$ in $H_{k} \setminus H_{k-1}$. By linear independence, we have $g = (u_k - \lambda_k)q_k  + \sum_{i\leq k-1} q_{i} v_i$ with $q_k \neq 0$. Set $g_0 = g + q_k \lambda_k$; this is in $H_{k-1} + u_k \Q$.

Suppose $s_{-} (g) = \alpha$, so that $g \geq \alpha$ (as functions on $K$), and thus $g_0 \geq \alpha + \lambda_k q_k$. There exists pure $\tau_0$ \st $\tau_0 (g) = s_- (g)$; obviously, 
$\tau_0 (g_0) = \alpha + \lambda_k q_k$. Thus $s_-(g_0) = \alpha + \lambda_k q_k$, so the latter number is in $s(H_{k-1} + u_k \Q)$. Hence $s_-(g) = \alpha = -\lambda_k q_k + (\alpha + \lambda_k q_k)$, so belongs to the coset, $-\lambda_k q_k+ s(H_{k-1} + \Q u_i)$, in $\R/s(H_{k-1} + \Q u_k)$; this is nonzero since $q_k$ and $\lambda_k$ are  not zero and  $\lambda_k \Q \cap s(H_{k-1} + \Q u_k) = \brcs{0}$. 

The same computation with a couple of inequalities reversed  shows that $s_+ (g)\in -\lambda_k q_k+ s(H_{k-1} + \Q u_k)$. This completes the inductive process.

Now set $G$ to be the rational span of $\brcs{v_i}_{0}^{\infty}$; obviously, this is an increasing union of the $H_k$. Since $G$ is a rational vector space, its closure (within $A$) is a real vector space, and therefore contains $\R G$, which in turn contains the real span of $\brcs{u_i}$, and is thus dense. Hence $G$ is dense in $A$. Impose the strict ordering on $G$, so that $G$ is a simple dimension group (we finally use the fact that $K$ is a Choquet simplex). Now we verify that $G$ satisfies the condition of Lemma 1. 

Select nonzero $g$ in $G$. Since $G$ is the union of the ascending chain of subgroups $\brcs{H_i}$, there exists $k \geq 0$ \st $g$ belongs to $H_k \setminus H_{k-1}$ (define $H_{-1} = \brcs{0}$). From the computation above, $s_- (g)$  belongs to a nontrivial coset of a rational subgroup of the reals, hence cannot be zero. Thus $G$ is a simple archimedean dimension group,  a dense subgroup of $\Aff K$, and it is immediate that its normalized trace space is $K$. \qed

\comment 
\subtitle References

\long\def\Reff[#1] #2, #3, #4\par{\vskip 1pt \item{[#1]} #2, {\it #3,}
#4\par} {\parindent = 2.5em 

\Reff [EHS] Edward G Effros{, David E Handelman,} and Chao-Liang Shen,
Dimension Groups and Their Affine Representations, 
American Journal of Mathematics 102 (1980) 385--407.

\Reff [G] K Goodearl, Partially ordered abelian groups with interpolation, Mathematical series and monographs \#20, American Mathematical Society (1986).

\Reff [GH] K Goodearl and D
Handelman, Tensor products of dimension groups and K${}_0$ of unit regular rings, Canad J Math 38 (1986) 633--658.

}  

\endcomment  
\vskip 6pt \noindent Mathematics Department, University of Ottawa, Ottawa
ON K1N 6N5, Canada; dehsg\@uottawa.ca

\end

%% file: generic_macros




\font\rm=cmr10 \rm

\font\bf=cmb10
\font\Rm=cmr9 at 11pt
\rm
\font\it=cmsl9 at 10pt
 at 7pt

\font\Rrm=cmr17 at 16pt
   \font\Rm=cmr12 at 11.5pt

\long\def\Pf{\par\noindent {\it Proof.} }
\def\({\left(}
\def\){\right)}
\def\st{such that }
\def\qed{\hfill$\bullet$\vskip 4pt}

\def\brcs#1{\left\{ #1\right\}}

\def\wrt{with respect to }
\def\:{\,:}

\def\R{\text{\bf R}}
\def\N{\text{\bf N}}
\def\Z{\text{\bf Z}}
\def\Q{\text{\bf Q}}

\def\Arrow #1;#2.{#1\:#2 \to }

\def\Set#1#2{\brcs{#1 \left|\vphantom{#1 #2} \right.#2}}



\def\Rrr#1,#2{{\Cal J}_{#1,#2}}
\def\slfrac#1#2{{\raise -.07 ex\hbox{$^{#1}$}}\!/\raise .35 ex \hbox{${}_{#2}$}}
\def\ssf #1/#2{\slfrac {#1}{#2}}

\def\pd #1,#2.{\frac {\partial #1}{\partial #2}}

   \long\def\Lem
#1.#2\par{\vskip4pt{\baselineskip=13pt\font\it=cmsl12 at
11.5pt\Rm
   \noindent {\rm \uppercase{#1}} #2\vskip3pt

   }}

\long\def\Title #1\par {\noindent{\Rrm #1}\vskip 9pt}

 \long\def\SubT #1.{\noindent {\it #1\/} } 
 
 \long\def\SecT
#1\par{\vskip 3pt \noindent {\bf #1}\vglue1pt
   \noindent}

\long\def\subtitle #1.{\vskip 2pt \noindent {\it #1}}

\long\def\Rmk#1\par{\vskip 1pt \noindent {\it
Remark.} #1\vskip2pt}


%% file: papermacros_entirefunctions
\scrollmode\NoBlackBoxes
\magnification=1100
\long\def\Abstract #1\par%
{\vskip .2 true cm{\leftskip 1 true in \rightskip 1 true in \font\rm=cmr8 \rm
\baselineskip=1pt \font\it=cmsl8 \font\bf=cmb8
\parindent=0em {\bf Abstract} #1

}}
\comment
\font\rm=Times at 10pt

\font\bf=TimesB
\font\Rm=Times at 11pt
\rm
\font\it=TimesI at 10pt
\endcomment

\long\def\Pf{\par\noindent {\it Proof.} }
\def\({\left(}
\def\){\right)}
\def\st{such that }
\def\qed{\hfill$\bullet$\vskip 4pt}

\def\brcs#1{\left\{ #1\right\}}
\def\Set#1#2{\brcs{#1 \left|\vphantom{#1 #2} \right.#2}}

\def\wrt{with respect to }
\def\:{\,:}
\def\Arrow #1;#2.{#1\:#2 \to }


\def\R{\text{\bf R}}
\def\N{\text{\bf N}}
\def\Z{\text{\bf Z}}
\def\Q{\text{\bf Q}}
 
\def\Rrr#1,#2{{\Cal J}_{#1,#2}}

\def\slfrac#1#2{{\raise -.07 ex\hbox{$^{#1}$}}\!/\raise .35 ex \hbox{${}_{#2}$}}
\def\ssf #1/#2{\slfrac {#1}{#2}}

\def\pd #1,#2.{\frac {\partial #1}{\partial #2}}


   \long\def\Title #1\par {\noindent{\Rrm #1}\vskip 9pt}
 \long\def\SubT #1.{\noindent {\it #1\/} }   \long\def\SecT
#1\par{\vskip 3pt \noindent {\bf #1}\vglue1pt
   \noindent}
\long\def\subtitle #1.{\vskip 2pt \noindent {\it #1}}

\long\def\Rmk#1\par{\vskip 1pt \noindent {\it
Remark.} #1\vskip2pt}

